\documentclass[12pt,leqno,fleqn,epsfig]{article}
\usepackage{amssymb, epsfig, amsmath, amsthm}
\usepackage{mathrsfs}       
     \usepackage{esint}   
     
\newcommand{\R}{\mathbb{R}}

\newcommand{\N}{\mathbb{N}}

\newcommand{\mes}{\operatorname{\rm meas}}    
\newcommand{\esssup}{\operatorname*{ess\,sup}}

\newcommand{\const}{\operatorname*{const}}

\newcommand{\bb}{\begin{equation}}
\newcommand{\ee}{\end{equation}}
\newcommand{\bq}{\begin{eqnarray}}
\newcommand{\eq}{\end{eqnarray}}
\newcommand{\bqn}{\begin{eqnarray*}}
\newcommand{\eqn}{\end{eqnarray*}}

\newcommand{\var}{\varepsilon}

\newcommand{\intl}{\int\limits}

\newcommand{\Beweisende}{\rule{0.2cm}{0.2cm}}

\newcommand{\D}{\displaystyle}

\newcommand{\intmw}{{\int\hspace{-830000sp}-\!\!}}

\newcounter{secnum}

\newtheorem{thm}{Theorem}[section]
\newtheorem{cor}[thm]{Corollary}
\newtheorem{lem}[thm]{Lemma}

\newtheorem{prop}[thm]{Proposition}

\theoremstyle{definition}

\newtheorem{defin}[thm]{Definition}
\newtheorem{rem}[thm]{Remark}

\title{On the Liouville type theorems for  self-similar  solutions to the Navier-Stokes equations} 
 
\author{Dongho Chae$^*$  and J\"{o}rg Wolf $^\dagger$\\
\ \\
 $*$Department of Mathematics\\
Chung-Ang University\\
 Seoul 156-756, Republic of Korea\\
 e-mail: dchae@cau.ac.kr\\
and \\
$\dagger$Department of Mathematics\\
Humboldt University Berlin\\
Unter den Linden 6, 10099 Berlin, Germany\\
e-mail: jwolf@math.hu-berlin.de}
\date{}
\begin{document}
\maketitle
\begin{abstract}
We prove Liouville type theorems for the self-similar solutions to the Navier-Stokes equations.  One of our results generalizes the previous ones by Ne\v{c}as-R\.{u}\v{z}i\v{c}ka-\v{S}verak and Tsai. 
Using the Liouville type theorem we also remove a scenario of asymtotically self-similar blow-up for the Navier-Stokes equations with the profile belonging to $L^{p, \infty} (\Bbb R^3)$ with $p> \frac{3}{2}$.
  \noindent{\bf
keywords:} Navier-Stokes equations, self-similar solution, Liouville's theorem

\end{abstract}

\section{Introduction}
\label{sec:-1}
\setcounter{secnum}{\value{section} \setcounter{equation}{0}
\renewcommand{\theequation}{\mbox{\arabic{secnum}.\arabic{equation}}}}

We consider the  Navier-Stokes equation in the space time cylinder  $ \R^{3}\times (-T,0)$
\begin{equation}
\partial _t u  + (u \cdot \nabla ) u -\Delta  u = - \nabla p,\qquad  \nabla \cdot  u =0,
\label{1.1}
\end{equation}
where $u=(u_1 (x,t), u_2 (x,t), u_3 (x,t))$,  $(x,t)\in  \Bbb R^3\times (-T, 0)$.  The aim of the present paper is to 
exclude a possible self similar blow up at the point $ (0,0)$ under more general assumptions than in \cite{tsa}.  More precisely,     
 we assume that $ u$ and $ p$ respectively  are  given by a self similar  profiles $ U: \R^{3} \rightarrow \R^{3} $ 
 and $ P: \R^{3} \rightarrow \R$ such that 
 \begin{align}
 u(x,t) &= \frac{1}{ \,\sqrt{-2a t}} U \Big(\frac{x}{ \,\sqrt{-2a t}}\Big),
 \label{1.2}
  \\
 p(x,t) &= \frac{1}{ \,-2a t} P \Big(\frac{x}{ \,\sqrt{-2a t}}\Big),\quad  (x,t) \in  \R^{3}\times (-T,0), 
 \label{1.3}
  \end{align}
 where $ a $ is a positive constant. 
 Then $ (U, P)$ solves the following system  proposed by Leray (cf. \cite{ler}). 
\begin{equation}
-\Delta U  + (U \cdot \nabla ) U + a y\cdot \nabla U + a U = - \nabla P,\qquad  \nabla \cdot  U =0\quad  \text{ in}
\quad  \R^{3}. 
\label{1.4}
\end{equation}
It is already known that if $ U \in L^p(\R^{3})$ for some $ p \in [3, +\infty]$,  then $U=0$ for $p\in [3, \infty)$, while $ U = \const$ for $p=\infty$. The case $ p=3$ 
is proved in \cite{nec},  while the case $ p>3$ has been proved by Tsai in \cite{tsa}. In fact Tsai proved a more general result, namely that 
$ U = 0$ if $ u $ satisfies a the local energy bound 
\begin{equation}
\sup_{t\in (-t_0, 0)} \intl_{B}  |  u(t)|^2 dx  +  \intl_{-t_0}^{0} \intl_{B}  | \nabla  u|^2 dx dt < +\infty
\label{1.5}
\end{equation}
for some ball $ B \subset \R^{3}$ and some $ t_0 >0$. 

\hspace{0.5cm}
We extend the results mentioned above in different directions. Our first main result is the following 

\begin{thm}
\label{thm1.1}
Let $ (U, P) \in C^{\infty}(\R^{3})^3 \times C^{\infty}(\R^{3})$  be a solution to \eqref{1.4}, and $\Omega=\nabla \times U$.  Suppose that for some 
$ q>0$
\begin{align}
& \| U\|_{ L^q(B_1(y_0))}+ \|  \Omega\|_{ L^2(B_1(y_0))} = o(| y_0|^{ \frac{1}{2}}) \quad \mbox{as}\quad |y_0|\to +\infty.
\label{1.6}
\end{align}
Then,  $ U$ is a constant function. 
\end{thm}  

\hspace{0.5cm}
Below we remove the condition on $ \Omega $, and instead  we  restrict the range of $q$ so that $ q > \frac{3}{2}$.  
Our second main result is the following 
\begin{thm}
\label{thm1.2}
Let $ (U, P) \in C^{\infty}(\R^{3})^3 \times C^{\infty}(\R^{3})$  be a solution to \eqref{1.4}.  Suppose that for some $ \frac{3}{2} <q < +\infty$ and $ \alpha >0$
\begin{align}
\intl_{B_1(y_0)\cap \{ | U|>\alpha \}} | U|^q dx \to 0 \quad  \text{ as}\quad  | y_0| \rightarrow +\infty. 
\label{1.7}
\end{align}
Then $ U$ is a constant function. 
\end{thm}

\begin{rem}
\label{rem1.3} If $U\in L^\infty (\Bbb R^3)$, then (\ref{1.7}) is obviously satisfied with the choice of $\alpha=\|U\|_{L^\infty} +1$.
In general, if 
$ U\in L^{p.\infty}(\R^{3})$ for $ p>q$ implies \eqref{1.7}.  Indeed,  for   $\frac{q}{p}<\theta < 1$ we have 
\begin{align*}
\intl_{B_1(y_0)\cap \{ | U| > \alpha \}} | U|^q dx &= q \intl_{\alpha }^{\infty} \sigma ^{ q-1} \mes \{ B_1(y_0)\cap | U| >\sigma\}  d \sigma 
\\
&\le q \| U\|_{ L^{ p,\infty}}^{p \theta } \mes\{ B_1(y_0)\cap | U| > \alpha \}^{ 1-\theta }  \intl_{\alpha }^{\infty} \sigma ^{ q-p\theta -1} d\sigma 
\\
&= q\| U\|_{ L^{ p,\infty}}^{ p\theta } \mes\{ B_1(y_0)\cap | U| > \alpha \}^{ 1-\theta }\frac{\alpha^{ q-p\theta }}{p\theta -q}\to 0
\end{align*}
as $|y_0|\to +\infty$.
Thus,  Theorem\,\ref{thm1.2}  leads to the following Corollary 

\begin{cor}
\label{cor1.4}
Let $ (U, P) \in C^{\infty}(\R^{3})^3 \times C^{\infty}(\R^{3})$  be a solution to \eqref{1.4}.  Suppose that for some 
$ \frac{3}{2} < p < +\infty$ 
\begin{align}
U\in L^{p, \infty}(\R^{3}). 
\label{1.8}
\end{align}
\end{cor}

\end{rem}
The above corollary shows clearly that Theorem\, \ref{thm1.2} improves the previous results of \cite{nec, tsa}.
As an application  the above result one can remove a scenario of asymptotically self-similar blow-up  with a profile given by (\ref{1.8}) as follows, which could viewed as an improvement of the corresponding result in \cite{cha}.
\begin{thm}
\label{thm1.3}
Let $u\in C^2 (\Bbb R^3 \times (0, t_*))$ be a solution to (\ref{1.1}).  Suppose there exists $U$ satisfying (\ref{1.8}) with  $ \frac{3}{2}< p < +\infty$,  and $ q \ge 2$ such  that 
\begin{equation}
\label{1.9}
\lim_{t\nearrow t_*}(t_*-t)^{\frac{q-3}{2q} } \sup_{t<\tau <t_*}\left \|u (\cdot , \tau)- \frac{1}{\sqrt{2a(t_*-\tau)}} U\left( \frac{\cdot -x_*}{\sqrt{ 2a(t_*-\tau )}} \right) \right\|_{L^q (B_{r\sqrt{t_*-t} } (x_*))} =0
\end{equation}
for all $r>0$. Then, $U=0$, and $z_*=(x_*,t_*)$ is not a  blow-up  point.
\end{thm}

\section{Local $ L^\infty$ estimate for local suitable weak solutions to the Navier-Stokes equation without pressure}
\label{sec:-4}
\setcounter{secnum}{\value{section} \setcounter{equation}{0}
\renewcommand{\theequation}{\mbox{\arabic{secnum}.\arabic{equation}}}}

The aim of the present section is to provide a local  $ L^\infty$ bound for local suitable weak solutions to the Navier-Stokes 
equations without pressure.  

\hspace{0.5cm}
First, let us recall the notion  of the local pressure projection 
$ E^{ \ast}_{G}: W^{-1,\, s}_{ G} \rightarrow  W^{-1,\, s}(G)$ for a given bounded $ C^2$-domain $ G \subset \R^{n}$, 
$ n\in \R^{n}$, introduced in \cite{wol3}.  
Appealing to the $ L^p$-theory of the steady Stokes system (cf.  \cite{gal}), 
for any  $ F\in W^{-1,\, s}_{ G}$ there exists a unique pair $ (v, p)\in W^{1,\, s}_{ 0, \sigma }(G)\times L^s_0(G)$ 
 which solves in the weak sense  the steady Stokes system 
 \begin{align*}
 \nabla \cdot v &=0\quad  \text{ in}\quad  G,\quad  -\Delta v + \nabla p = F\quad  \text{ in}\quad  G,
 \\
&\qquad \qquad \qquad  v=0   \text{ on}\quad  \partial G.
 \end{align*}
Then we set $E^{ \ast}_G(F):=  \nabla p $, where $ \nabla p $ denotes the gradient functional in $ W^{-1,\, s}(G)$  defined by 
\[
 \langle \nabla p, \varphi \rangle = \intl_{G} p \nabla \cdot \varphi dx, \quad  \varphi \in W^{1,\, s'}_0(G). 
\]
Here we have denoted by $ L^s_0(G)$ the space of all $ f\in L^s(G)$ with $ \intl_{G} fdx=0 $.

\begin{rem}
\label{rem4.1}
1. The operator $ E^{ \ast}_G$  
is bounded from $ W^{-1,\, s}(G)$ into itself with  $ E^{ \ast}_G(\nabla p)=\nabla p$ for all $ p\in L^s_{ 0}(G)$. The norm of $ E^{ \ast}_G $ depends only on $ s$ and the geometric properties of $ G$, and independent on $ G$,  if $ G$ is a ball or an annulus,  which is due to the scaling properties of the Stokes equation.   

\hspace{0.5cm}
2. In case $ F\in L^s(G)$ using the canonical embedding  $ L^s(G) \hookrightarrow  W^{-1,\, s}(G)$ and  the  elliptic regularity 
 we get $ E^{ \ast}_G(F)= \nabla p \in L^s(G)$ together with the estimate 
\begin{equation}
\| \nabla p\|_{s, G} \le c \| F\|_{ s, G}, 
\label{4.1}
\end{equation}
where the  constant in \eqref{4.1}  
depends only on $ s$ and $ G$. In case  $ G$  is a ball or an annulus this constant depends only on $ s$ 
 (cf. \cite{gal} for more details). Accordingly the restriction of  $ E^{ \ast}_G$ to the Lebesgue space $ L^s(G)$ defines  
a projection in $L^s(G)$. This projection will be denoted still by $ E^{ \ast}_G$.  \\
\ \\
Below for a class of vector fields $X$ we denote
by $X_{\sigma}$ the set of $u\in X$  such that $\nabla \cdot u=0$ in the sense of distribution.
\end{rem}

\hspace{0.5cm}
By using the projection $ E^{ \ast}_G$, we introduce  the following notion of local suitable weak solution to the 
Navier-Stokes equations

\begin{defin}[Local suitable weak  solution] 
  \label{def1.1}
Let $ Q= \R^{3}\times (-T,0)$. A vector function
 $ u\in L^2_{loc }(Q)$ is called a {\it local  suitable weak  solution to \eqref{1.1}}, if 
  \begin{itemize}
  \item[1.]   $ u\in L^{ \infty}_{ loc}(-T,0; L^2_{ loc}(\R^{3}))\cap L^2(-T, 0; W^{1,\, 2}_{ loc, \sigma } (\R^{3}))$.
  \item[2.]  $u $ is a distributional solution to \eqref{1.1}, i.\,e. for every $ \varphi \in C^{\infty}_{\rm c}(Q)$  with $ \nabla \cdot u=0$
  \begin{equation}
  \intl  \hspace*{-0.2cm}\intl  _{ \hspace*{-0.5cm} Q}- u\cdot \frac{\partial \varphi }{\partial t} - u \otimes u : \nabla \varphi + \nabla u : \nabla \varphi dxdt =0. 
 \label{4.2}
 \end{equation}
 \item[3.] For every ball $ B \subset \R^{3}$ the following local energy inequality {\em without pressure}  holds  for  every nonnegative $ \phi \in C^{\infty}_{\rm c}(B\times (0,+\infty))$, and for almost every $ t\in (-T,0)$ 
 \begin{align}
 &   \frac{1}{2} \intl_{B} | v_B(t)|^2 \phi dx +    \intl_{-T}^{t}\intl_{B} | \nabla v_{B}|^2  \phi dx    ds
\cr
 & \quad  \le  \frac{1}{2}  \intl_{-T}^{t} \intl_{B} | v_{B}|^2 \Big(\Delta + \frac{\partial }{\partial t}\Big) \phi  + 
 | v_B|^2 u\cdot \nabla \phi )  dx    ds
 \cr
 & \qquad +  
 \intl _{-T}^{t}\intl_{B} (u \otimes v_{B}) :\nabla ^2 p_{h, B} \phi dx dt   +\intl_{-T}^{t} \intl_{G} p_{ 1,B } v_{B} \cdot \nabla \phi dxds
 \cr
 &\qquad \qquad \qquad +  
 \intl_{-T}^{t}\intl_{G} p_{2,B } v_{B}\cdot \nabla \phi dxds,
 \label{4.3}
\end{align} 
 where $ v_{B}= u + \nabla p_{ h, B }$, and 
\begin{align*}
\nabla p_{ h,B } &= -E^{ \ast}_{B} (u),
\\
\nabla p_{ 1,B } &= -E^{ \ast}_{ B} ((u\cdot \nabla ) u), \quad \nabla p_{2,B } = E^{ \ast}_{B} (\Delta u). 
\end{align*} 

  \end{itemize} 

 \end{defin}
 
 \begin{rem}
 1. Note that due to $ \nabla \cdot u=0$ the pressure $ p_{ h, B}$ is harmonic, and thus smooth in $ x$. 
 Furthermore, as it has been proved in \cite{wol3} the pressure gradient  $ \nabla p_{ h,B}$ is continuous 
 in $ B\times (-T,0)$. 
 
 \hspace{0.5cm}
 2. The notion of local suitable weak solutions to the Navier-Stokes equations satisfying the local energy inequality 
 \eqref{4.3}  has been  introduced in \cite{wol2}.  As it has been shown there such solutions enjoy the same partial regularity as the standard suitable weak solution 
as proved in the paper by  Caffarelli-Kohn-Nirenberg\cite{caf}. Furthermore, the following $ \var $-regularity criterion has been proved for solution satisfying 
 \eqref{4.3}:
 
\begin{itemize}
  \item[] {\it There exists and absolute number $ \var >0$ such that if for any $ Q_r = B_r(x_0)\times (t_0, t_0- r^2)$ it holds 
  \[
  r^{ -2}  \intl  \hspace*{-0.2cm}\intl  _{ \hspace*{-0.5cm} Q_r}| u|^3 dxdt \le 
    \var ^3 \quad  \Longrightarrow \quad    
    u\in L^\infty(Q_{ r/2})
\]}
\end{itemize} 
(cf. also \cite{wol2}).
\end{rem}
\hspace{0.5cm}

\hspace{0.5cm}
Before turning to the statement of this result we will fix the notations used throughout this section  
For  $ z_0 =(x_0, t_0)\in Q$ and $ 0<r< \,\sqrt{-t_0}$ we define the parabolic 
cylinders 
\[
Q_r = Q_r(z_0) = B_r(x_0)\times (t_0-r^2, t_0),\quad I_r= I_r(t_0)=(t_0-r^2, t_0). 
\]
By $ V^2_\sigma (Q_r)$ we denote the space $L^\infty(I_r; L^{ 2}(B_r))\cap  L^2(I_r; W^{1,\, 2}_\sigma (B_r ))$. 
Furthermore for $ u\in V^2(Q_r(z_0))$ we set
\begin{align*}
\hspace*{-1cm}A_q(r, z_0) &= \bigg(r^{q-3}\esssup_{t\in I_r(t_0)} \intl_{B_r(x_0)} | u(x,t)|^q dx\bigg)^{ \frac{1}{q}},\quad  G(r, z_0)= r^{ -1} \intl_{Q_r(z_0)} | \nabla u|^2dxdt, 
\\
&\qquad E_q( r, z_0) = r^{ -\frac{3+q}{3q}}   
\bigg(\intl_{I_r(t_0)} \bigg(  \intl_{B_r(x_0)} | u|^q dx \bigg)^{ \frac{3}{2q-3}}dt\bigg)^{ \frac{2q-3}{3q}}. 
\end{align*}

\begin{rem}
\label{rem2.4}
According to Lemma\,4.1\cite{wol2} the following Caccioppoli-type inequality holds true 
\begin{equation}
G \Big(\frac{r}{2}, z_0\Big)\le C \Big(E_3(r, z_0)^2 + E_3(r, z_0)^{ 3}\Big),
\label{4.4}
\end{equation}
where $ C>0$ denotes an absolute constant. 
\end{rem}

\hspace{0.5cm}
Our main result of this section  is the following $ \var $-regularity criterion  

\begin{thm}
\label{thm4.5}
Let $ u\in V^2_\sigma (Q)$ be a local suitable weak solution to \eqref{1.1}. 
Let $ \frac{3}{2} < q \le 3$. There exist two positive  constants $ \var _{q}$ and $ C_{q}$, both depending on $ q$ only,  such that if 
for $ Q_r(z_0) \subset Q$, $ z_0=(x_0, t_0)$,  the condition 
\begin{equation}
A_q(r, z_0)\le \var _{q}. 
\label{4.5}
\end{equation}
implies $ u\in  L^\infty (Q_{ \frac{r}{2}}(z_0))$, and it holds 
\begin{align}
 \esssup_{Q_{ \frac{r}{2}}(z_0)} | u|
\le C_{ q} r^{ -1}A_q(r, z_0).
\label{4.7}
\end{align}
\end{thm}

\hspace{0.5cm}
Before turning to the proof of Theorem\,\ref{thm4.5} we provide some lemmas, which will be used in our discussion below. 
We begin  with a Caccioppoli-type inequlities similar to \eqref{4.4}. 

\begin{lem}
\label{lem4.6}
Let $ u\in V^2_\sigma (Q_R)$ be local suitable weak solution to \eqref{1.1}. Then for every $ \frac{3}{2} < q \le 3$
\begin{equation}
E_3\Big(\frac{3}{4}R\Big)^2 + G \Big(\frac{3}{4}R\Big) \le 
C \Big(E_q(R)^2 + E_q(R)^{ \frac{3q}{2q-3}}\Big) \le C \Big(A_q(R)^2 + A_q(R)^{ \frac{3q}{2q-3}}\Big), 
\label{4.8}
\end{equation}
  where $ C>0$ denotes a constant depending only on $ q$. 
\end{lem}  

{\bf Proof}: Let $ 0< r< \rho  \le R$ be fixed. Set $ B= B_\rho $, and define $ v_B = u+ \nabla p_{ h,B}$, 
where $ \nabla p_{ h, B}= - E^{ \ast}_{ B}(u)$.  Let $ \phi $ denote a suitable cut off function for $ Q_r \subset Q_\rho $. 
As it has been proved in  \cite{wol2} (cf. estimate (4.4) therein), applying H\"older's inequality, the following inequality holds
\begin{align}
&  \| \phi v_B\|^2_{ L^\infty(I_\rho ; L^{ 2}(B_\rho ))} + \| \phi\nabla v_B \|^2_{ L^2(Q_\rho )}
\cr
&\qquad \le c \rho (\rho -r)^{ -2}\| u\|_{ L^3(I_\rho ; L^{ \frac{9}{4}}(B_\rho ))}^2 + (\rho -r)^{ -1}\| u\|_{ L^3(Q_\rho )}^3 
+ \frac{1}{4} \| \nabla u \|^2_{ L^2(Q_\rho )} 
\cr
&\qquad \le c \rho^{ \frac{5}{3}} (\rho -r)^{ -2}\| u\|_{ L^3(Q_\rho )}^2 + (\rho -r)^{ -1}\| u\|_{ L^3(Q_\rho )}^3 
+ \frac{1}{4} \| \nabla u \|^2_{ L^2(Q_\rho )} . 
\label{4.9a}
\end{align}
By means of Sobolev's inequality together with H\"older's inequality, \eqref{4.1}  and  \eqref{4.9a} 
\begin{align}
\| \phi v_B\|^2_{ L^3(I_\rho ; L^{ \frac{18}{5}}(B_\rho ))}
&\le \|\phi v_B \|_{L^\infty(I_\rho; L^2 (B_\rho ))} ^{\frac43} \|\phi v_B \|_{L^2 (I_\rho; L^6(B_\rho))} ^{\frac23}\cr
 &\le \|\phi v_B \|_{L^\infty(I_\rho; L^2 (B_\rho ))} ^{\frac43} (\|v_B \cdot \nabla \phi \|_{L^2(Q_\rho )} 
+\|\phi \nabla v_B\|_{L^2 (Q_\rho )} ) ^{\frac23} \cr
&\le
 \| \phi v_B\|^{ \frac{4}{3}}_{ L^\infty(I_\rho ; L^2(B_\rho ))} 
 \Big((\rho -r)^{ -1} \| v_B \|_{ L^2(Q_\rho )}
 + \| \phi\nabla v_B \|_{ L^2(Q_\rho )}\Big)^2 \cr
 &\le  c \| \phi v_B\|^{2}_{ L^\infty(I_\rho ; L^2(B_\rho ))}  +(\rho -r)^{ -2} \| v_B \|_{ L^2(Q_\rho )}^2 +\frac{1}{16}
  \| \phi\nabla v_B \|_{ L^2(Q_\rho )}^2\cr
&\le c \rho^{ \frac{5}{3}} (\rho -r)^{ -2}\| u\|_{ L^3(Q_\rho )}^2 + (\rho -r)^{ -1}\| u\|_{ L^3(Q_\rho )}^3 
+ \frac{1}{16} \| \nabla u \|^2_{ L^2(Q_\rho )} . 
\label{4.9}
\end{align}
We recall the following Caccioppoli inequality for a harmonic function
$$
\int_{B_\rho}\phi^2 |\nabla h|^2 dx \leq  \max_{B_\rho} |\Delta \phi|  \int_{B_\rho} |h|^2 dx,
$$
which will be repeatedly used below.
The proof is immediate from the formula 
$-\Delta h^2  +2 |\nabla h|^2 =0,
$
by multiplying $\phi$, integrating over $B_\rho$,  and then using  integration by part.
Recalling that $ p_{ h,B}$ is harmonic, by using \eqref{4.1} with $ s=3$  we get first
\begin{align}
\|\phi \nabla p_{h, B} \|_{ L^{ \frac{18}{5}}(B_\rho )} ^3
&=c\left( \int_{B_\rho} \phi ^{\frac{18}{5}} |\nabla p_{h,B} |^{\frac{18}{5}} dx \right) ^{\frac56}\leq c\left( \int_{B_\rho} \phi ^{6} |\nabla p_{h,B} |^{6} dx\right) ^{\frac12}\rho\cr
&\leq c(\rho-r)^{-3} \left(\int_{B_\rho}  |\nabla p_{h,B} |^{2}dx \right) ^{\frac32}\rho
+c \left( \int_{B_\rho} \phi ^{2} |\nabla^2 p_{h,B} |^{2}dx \right) ^{\frac32}\rho\cr
& \leq c(\rho-r)^{-3} \left(\int_{B_\rho}  |\nabla p_{h,B} |^{2} dx\right) ^{\frac32}\rho,
\end{align}
from which, integrating it over $I_\rho$,  we obtain
$$
\| \phi \nabla p_{ h,B}\|^2_{ L^3(I_\rho ; L^{ \frac{18}{5}}(B_\rho ))}
\le \rho ^{ \frac{5}{3}} (\rho -r)^{ -2}
c\| u\|^2_{ L^3(Q_\rho )}.
$$
Using this estimate, we have
\begin{align}
 \| \phi u\|^2_{ L^3(I_\rho ; L^{ \frac{18}{5}}(B_\rho ))} &\le 
 2\| \phi v_B\|^2_{ L^3(I_\rho ; L^{ \frac{18}{5}}(B_\rho ))} + 2\| \phi \nabla p_{ h,B}\|^2_{ L^3(I_\rho ; L^{ \frac{18}{5}}(B_\rho ))}
 \cr
&\le  2\| \phi v_B\|^2_{ L^3(I_\rho ; L^{ \frac{18}{5}}(B_\rho ))} + c\rho ^{ \frac{5}{3}} (\rho -r)^{ -2}
c\| u\|^2_{ L^3(Q_\rho )}.
\label{4.9b}
\end{align}
Combining \eqref{4.9} with \eqref{4.9b}, we arrive at 
 \begin{align}
\| u\|^2_{ L^3(I_r ; L^{ \frac{18}{5}}(B_r ))}
\le c \rho^{ \frac{5}{3}} (\rho -r)^{ -2}\| u\|_{ L^3(Q_\rho )}^2 + (\rho -r)^{ -1}\| u\|_{ L^3(Q_\rho )}^3 
+ \frac{1}{8} \| \nabla u \|^2_{ L^2(Q_\rho )} . 
\label{4.9c}
\end{align} 
Once more using the fact that $ p_{ h,B}$ is harmonic applying integration by parts, Caccioppoli type inequality together \eqref{4.1},   we evaluate for almost all 
$t \in I_\rho $
\begin{align*}
\|\phi(t)  \nabla u(t)\|_{ L^2(B_\rho )}^2 &= 
\intl_{B_\rho } | \nabla v_B(t)|^2 \phi ^2(t) dx + 
\intl_{B_\rho } (\nabla v_B(t)+ \nabla u(t)) : (\nabla v_B(t)- \nabla u(t))\phi ^2(t) dx 
\\
&=\intl_{B_\rho } | \nabla v_B(t)|^2 \phi ^2(t) dx + 
\intl_{B_\rho } (\nabla v_B(t)+ \nabla u(t)): \nabla ^2 p_{ h,B} \phi ^2(t) dx 
\\
&= \intl_{B_\rho } | \nabla v_B(t)|^2 \phi ^2(t) dx -
\intl_{B_\rho } (v_B(t)+u(t)) \otimes \nabla \phi ^2(t): \nabla ^2 p_{ h,B}dx 
\\
& \le \|\phi (t) \nabla v_B(t) \|_{ L^2(B_\rho )}^2 + c (\rho -r)^{ -2} \| u(t)\|^2_{ L^2(B_\rho )}. 
\end{align*}
Integration of both side of the above inequality together with H\"older's inequality gives 
\begin{equation}
\|\phi  \nabla u\|_{ L^2(Q_\rho )}^2 \le \|\phi  \nabla v_B\|_{ L^2(Q_\rho )}^2
+ c\rho^{ \frac{5}{3}} (\rho -r)^{ -2}\| u\|^2_{ L^3(Q_\rho )}. 
\label{4.9d}
\end{equation}
Combining \eqref{4.9a} with \eqref{4.9d} we are led to
\begin{align}
  \| \nabla u \|^2_{ L^2(Q_r )}
 \le c \rho^{ \frac{5}{3}} (\rho -r)^{ -2}\| u\|_{ L^3(Q_\rho )}^2 + (\rho -r)^{ -1}\| u\|_{ L^3(Q_\rho )}^3 + \frac{1}{4} \| \nabla u \|^2_{ L^2(Q_\rho )} . 
\label{4.9e}
\end{align}
  Thus, adding \eqref{4.9c}  to \eqref{4.9e}, we obtain    
\begin{align}
&   \| u\|^2_{ L^3(I_r; L^{ \frac{18}{5}}(B_r ))} + \| \nabla u \|^2_{ L^2(Q_r )}
\cr
&\qquad \le c \rho^{ \frac{5}{3}} (\rho -r)^{ -2}\| u\|_{ L^3(Q_\rho )}^2 + (\rho -r)^{ -1}\| u\|_{ L^3(Q_\rho )}^3 
+ \frac{3}{8}\| \nabla u \|^2_{ L^2(Q_\rho )}. 
\label{4.10}
\end{align}  

Let $ t\in I_\rho $ be chosen  so that $ u(t)\in W^{1,\, 2}(B_\rho )$.  Applying H\"older's inequality together with Poincar\'{e}-Sobolev's inequality,  we see that 
\begin{align*}
\| u(t)\|^3_{ L^3(B_\rho )} &\le  \| u(t)\|_{ L^q(B_\rho )}^{ \frac{3q}{6-q}} \| u(t)\|_{ L^6(B_\rho )} ^{ 2\frac{9-3q}{6-q}}
\\
&\le  c\| u(t)\|_{ L^q(B_\rho )}^{ \frac{3q}{6-q}} \|\nabla  u(t)\|_{ L^2(B_\rho )} ^{ 2\frac{9-3q}{6-q}}
+  c\rho ^{ \frac{3q-9}{q}}\| u(t)\|_{ L^q(B_\rho )}^{3}. 
\end{align*}  
Integrating this inequality over $ I_\rho $, and applying H\"older's inequality, we are led to 
\begin{equation}
\| u\|_{ L^3(Q_\rho )}^3 \le c \| u\|_{ L^{ \frac{3q}{2q-3}}(I_\rho ; L^q(B_\rho ))}^{ \frac{3q}{6-q}}
 \| \nabla u\|_{ L^2(Q_\rho )}^{ 2\frac{9-3q}{6-q}} + 
 c\rho ^{ \frac{q-3}{q}}\| u\|_{ L^{ \frac{3q}{2q-3}}(I_\rho ; L^q(B_\rho ))}^{3}.
\label{4.11}
\end{equation}
We now estimate the right-hand side of \eqref{4.10} by the aid of \eqref{4.11},  and applying Young's inequality. This gives 
\begin{align}
&   \| u\|^2_{ L^3(I_r; L^{ \frac{18}{5}}(B_r ))} + \| \nabla u \|^2_{ L^2(Q_r )}
\cr
&\qquad \le c \rho^{ \frac{30-5q}{3q}} (\rho -r)^{ \frac{2q-12}{q}}
 \| u\|_{ L^{ \frac{3q}{2q-3}}(I_\rho ; L^q(B_\rho ))}^{2} + c(\rho -r)^{-  \frac{6-q}{2q-3}}
  \| u\|_{ L^{ \frac{3q}{2q-3}}(I_\rho ; L^q(B_\rho ))}^{ \frac{3q}{2q-3}} 
 \cr 
& \qquad \qquad +   c\rho ^{ \frac{q-3}{q}} (\rho -r)^{ -1} \| u\|_{ L^{ \frac{3q}{2q-3}}(I_\rho ; L^q(B_\rho ))}^{ 3}
+ \frac{1}{2}\| \nabla u \|^2_{ L^2(Q_\rho )}. 
\label{4.12}
\end{align}     
By using a standard iteration argument (e.g. see \cite{gia}) we deduce from \eqref{4.12} together with Young's inequality  that 
\begin{align}
&   \| u\|^2_{ L^3(I_{ \frac{3}{4}R}; L^{ \frac{18}{5}}(B_{ \frac{3}{4}R} ))} + \| \nabla u \|^2_{ L^2(Q_{ \frac{3}{4}R} )}
\cr
&\qquad \le c R^{ \frac{q-6}{3q}}  \| u\|_{ L^{ \frac{3q}{2q-3}}(I_R ; L^q(B_R))}^{2} + cR^{-  \frac{6-q}{2q-3}}
  \| u\|_{ L^{ \frac{3q}{2q-3}}(I_R ; L^q(B_R ))}^{ \frac{3q}{2q-3}}.  
 \label{4.13}
\end{align} 
Multiplying both sides of \eqref{4.13} by $ R^{ -1}$, and applying H\"older's inequality, we obtain the desired inequality \eqref{4.8}. 
 \hfill \Beweisende    

\vspace{0.5cm}  
\hspace{0.5cm}
We continue our discussion with some useful iteration lemmas. 
Let $ G \subset \R^{n}$ be a bounded $ C^2$-domain.  By $ A^s(G)$, $ 1< s < +\infty$,  we denote the image of 
$ W^{2,\, s}_0(G )$ under the Laplacian $ \Delta $, which is a closed subspace of $ L^s(G)$.  By $ B^s(G )$ 
we denote the complementary space, which contains all  $ p\in L^s(\Omega )$  being harmonic in $ G $ such that 
\begin{equation}
L^s(G) = A^s(G)+ B^s(G). 
\label{4.14}
\end{equation} 
By using the well-known Calder\'on-Zygmund inequality, and the elliptic regularity of the Bi-harmonic equation we get the following 

\begin{lem}
\label{lem4.7}
1. Let $ A \in L^s(G; \R^{n^2} )$. Then there exists a unique $ p_0 \in A^s(G )$ such that 
\begin{equation}
-\Delta p_0 = \partial _i \partial _j A_{ ij}\quad  \text{ in}\quad  G
\label{4.15}
\end{equation}
in the sense of distributions \,\footnotemark. 
\footnotetext{ Here \eqref{4.15} means $ - \intl_{G} p_0 \Delta \phi = \intl_{G} A_{ ij} \partial _{ i}\partial _j \phi $ for all $ \phi \in C^{\infty}_{\rm c}(G)$.  } 
In addition, it holds 
\begin{equation}
\| p_0\|_{s}  \lesssim  \| A\|_s.
\label{4.16}
\end{equation}

\hspace{0.5cm}
2. Let $ h \in L^s(G; \R^{n})$, $ 1 \le  s < n$. Then  there exists a unique $ p_0\in A^{ s^\ast}(G) 
\cap W^{1,\, s}(G)$ such that 
\[
-\Delta p_0 = \partial _i h_i\quad  \text{ in}\quad  G 
\]
in the sense of distributions, and the following estimate holds true
\begin{equation}
\| p_0\|_{ s^\ast}+ \|\nabla  p_0\|_{s}  \lesssim  \| h\|_s.
\label{4.18}
\end{equation}
The hidden constants in both \eqref{4.16} and  \eqref{4.18}  depend only on $ s, n$, and the geometric property of 
$G $. 
In case $ G$ equals a ball, these constants are independent of the radius.  

\end{lem}

\begin{lem}
\label{lem4.8}
Let $ f\in L^{ \frac{3}{2}}(Q_{ 1})$. Let $ 0<r_0 < 1$. Suppose, there exists 
$ 4 \le  \lambda \le 5$ and $ C>0$,  such that for all  $ z_0=(x_0, t_0) \in Q_{ \frac{1}{2}}$ and 
$ r_0 \le r \le \frac{1}{2}$
\begin{equation}
\intl_{Q_r(z_0)} |f- {\tilde f} _{ B_r(x_0)}|^{ \frac{3}{2}} \le   K_0 r^{ \lambda },
\label{4.18a}
\end{equation}
where $ {\tilde f} _{ B_r}(t) = {\D \intmw_{B_r} f(x,t) dx  }$. 
Let $ \nabla p = E_{ B_{ \frac{3}{4}}}^{ \ast}(\nabla \cdot f)$. Then for all $ z_0\in Q_{ \frac{1}{2}}$ and $ r_0 \le r \le  \frac{1}{4}$ it holds 
\begin{equation}
\intl_{Q_r(z_0)} | p-{\tilde p} _{ B_r(x_0)}|^{ \frac{3}{2}}  \le CK_0 r^{ 4 }.
\label{4.19}
\end{equation}
\end{lem}

{\bf Proof}:  Let $ z_0\in Q_{ \frac{1}{2}}$ and  $ r_0 \le r \le \frac{1}{8}$ be arbitrarily chosen, but fixed. Let $ 0< \theta < \frac{1}{2}$, 
specified below.  According to \eqref{4.14}  there exist unique  $ p_{ 0, r}(t) \in A^{ \frac{3}{2}}(B_r(x_0))$ and 
$ p_{ h,r}(t)\in B^{ \frac{3}{2}}(B_r(x_0))$ such that $ p(t)- {\tilde p} _{ B_r(x_0)}(t) = 
p_{ 0,r}(t) + p_{ h,r}(t)$. Noting that $ p(t)- {\tilde p}  _{ B_{ \theta r}(x_0)}(t) = 
 p(t)- {\tilde p}  _{ B_{r}(x_0)}(t) -  (p(t)- {\tilde p}  _{ B_{r}(x_0)}(t))_{ B_{ \theta r}(x_0)}$, it follows that 
 \begin{align*}
& \intl_{B_{ \theta r}(x_0)} | p(t)- {\tilde p} _{ B_{ \theta r}(x_0)}(t)|^{ \frac{3}{2}}  
 \\
&\qquad  \lesssim  
     \intl_{B_{ \theta r}(x_0)} | p_{ h,r}(t)- (p_{ h,r}(t))_{ B_{ \theta r}(x_0)}|^{ \frac{3}{2}}
     + \intl_{B_{ \theta r}(x_0)} | p_{ 0,r}(t)- (p_{ 0,r}(t))_{ B_{ \theta r}(x_0)}|^{ \frac{3}{2}}
\\
&\qquad  \lesssim   
\theta ^{ \frac{9}{2} }\intl_{B_r(x_0)} | p_{h,r}(t)|^{ \frac{3}{2}} + \intl_{B_r(x_0)} | f(t) - {\tilde f} _{ B_r(x_0)}(t)|^{ \frac{3}{2}}
\\    
& \qquad  \lesssim  \theta^{ \frac{9}{2}} \intl_{B_{ r}(x_0)} | p(t)- {\tilde p} _{ B_{ r}(x_0)}(t)|^{ \frac{3}{2}}  
+  \intl_{B_r(x_0)} | f(t) - {\tilde f} _{ B_r(x_0)}(t)|^{ \frac{3}{2}}. 
 \end{align*}
Integrating the above estimate over $ I_{ \theta r}(t_0)$, and observing the assumption \eqref{4.18},    we arrive at 
\begin{equation}
\intl_{Q_{ \theta r}(z_0)} | p- {\tilde p} _{ B_{ \theta r}(x_0)}|^{ \frac{3}{2}}\le 
 C_1\theta^{ \frac{9}{2}} \intl_{Q_{ r}(z_0)} | p- {\tilde p} _{ B_{ r}(x_0)}|^{ \frac{3}{2}}  + C_2 K_0 r^{4 }.
\label{4.20}
\end{equation}
By a standard iteration argument from \eqref{4.20} we deduce that 
\begin{equation}
\intl_{Q_r(z_0)} | p-{\tilde p} _{ B_r(x_0)}|^{ \frac{3}{2}}   \lesssim    r^4 \intl_{Q_{ \frac{1}{4}}(z_0)} | p|^{ \frac{3}{2}} + K_0 r^{ 4 }.
\label{4.21}
\end{equation}
Noting that by the definition of $ p$ having for almost every $ t\in I_{ \frac{1}{4}}(t_0)$
\[
\| p(t)\|_{ L^{ \frac{3}{2}}(B_{ \frac{1}{4}}(x_0))} \le c \| f(t)- {\tilde f} _{ B_{ \frac{3}{4}}}(t)\|_{ L^{ \frac{3}{2}}(B_{ \frac{3}{4}})},
\]
 the assertion \eqref{4.19} follows from \eqref{4.21} together with \eqref{4.18}.  \hfill \Beweisende

\vspace{0.5cm}  
\hspace{0.5cm}
We are in a position to prove the following iteration lemma, based on the idea of \cite{caf}.   

\begin{prop}
\label{prop4.9}
Let $ u\in V^2(Q_1(0,0))$ be a local suitable weak solution to the Navier-Stokes equations.  We define 
$ v= u+ \nabla p_h $, where $ \nabla p_h = - E^{ \ast}_{ B_{ \frac{3}{4}}}(u)$. 
There exist absolute positive numbers $ K_{ q}$ and $ \var_{ q}$ such that if 
\begin{equation}
A_q(1, 0) \le \var_{ q}
\label{4.22}
\end{equation}
then for all $ n \in \N, n \ge 2$, and for all $ z_0\in Q_{ \frac{1}{2}}(0,0)$ it holds 
\begin{align*}
& \hspace*{-2cm}(2.23)_n\qquad  \qquad \qquad \qquad \intmw_{Q_{ r_n}(z_0)} | v|^3 dxdt\le K_{ q}^3 A_q(1, 0)^3,
\end{align*}
where $ r_n= 2^{ -n}$, $n\in \N$. 
\end{prop}

 \setcounter{equation}{23}

{\bf Proof}:   From the definition of a local suitable weak solution the following local energy inequality holds 
true for every non negative $\phi \in C^{\infty}_{\rm c}\Big(B_{ \frac{3}{4}} \times \Big(- \frac{9}{16},0\Big]\Big)$, and for almost 
all $ t\in \Big(- \frac{9}{16},0\Big]$ 
\begin{align}
& \frac{1}{2}\intl | v(t)|^2 \phi (t)dx + \intl_{-r_3^2}^{t}\intl| \nabla v|^2  \phi dxds
\cr
&\qquad \le \frac{1}{2}\intl_{-r_3^2}^{t}\intl | v|^2  \Big(\frac{\partial \phi}{\partial t} + 
\Delta \phi\Big)dxds
+ \frac{1}{2}\intl_{-r_3^2}^{t}\intl | v|^2 v\cdot \nabla \phi - | v|^2 \nabla p_h\cdot \nabla \phi dxds
\cr
& \qquad\qquad   +\intl_{-r_3^2}^{t}\intl (v \otimes v - v \otimes \nabla p_h  : \nabla ^2 p_h) \phi dxds
+\intl_{-r_3^2}^{t}\intl p_1 v \cdot \nabla \phi dxds
\cr
&\qquad \qquad \qquad \qquad +\intl_{-r_3^2}^{t}\intl p_2 v \cdot \nabla \phi dxds.
\label{4.24}
\end{align}
where 
\begin{align*}
\nabla p_{ 1}&= E^{ \ast}_{ B_{  \frac{3}{4}}}(\Delta u), 
\quad  \nabla p_2= - E^{ \ast}_{ B_{  \frac{3}{4}}} (\nabla \cdot (u \otimes u)), 
\end{align*}
Note that by the definition of $ v$ it holds almost everywhere in $ Q_{ \frac{3}{4}}(0,0)$
\begin{equation}
u \otimes u =v \otimes v - v \otimes \nabla p_h - \nabla p_h \otimes v+ \nabla p_h \otimes \nabla p_h.
\label{4.25}
\end{equation}

\hspace{0.5cm}
\noindent{\it Proof of (2.23)$ _n$  by induction}: For $ n=2$ the inequality (2.23)$ _2$ follows immediately from Lemma\,\ref{lem4.6}.  

\hspace{0.5cm}
Let $ K_{ q}>1$ be a constant specified below.  
Assume $(2.23) _k$ is true for $ k=1, \ldots, n$, for some $ n\in \N$. This implies for all $ z_0\in Q_{ \frac{1}{2}}(0,0)$ and  $ r_{ n} \le r \le \frac{1}{2} $
\begin{equation}
\intmw_{Q_{ r}(z_0)} | v|^3 dxdt\le C K_{ q}^3 A_q(1, 0)^3.
\label{4.26a}
\end{equation}
Let $ r_{ n+1} \le r \le  r_3$ and $ z_0\in Q_{ \frac{1}{4}}(0,0)$ be arbitrarily chosen, but fixed. Using Cauchy-Schwarz's inequality, \eqref{4.26a}, and recalling that 
$p_h$ is harmonic, we get 
\begin{align}
\intmw_{Q_r(z_0)}  | v| ^{ \frac{3}{2}}| \nabla p_h|^{ \frac{3}{2}} dxdt& \le 
C K_{ q}^{ \frac{3}{2}}A_q(1, 0)^{ \frac{3}{2}} r^{ -1}
\Bigg[ \intl_{I_{ 1}}\bigg( \intl_{B_{ \frac{1}{4}}(x_0)} | \nabla p_h(t)| ^{q} dx\bigg) ^{ \frac{3}{q}} dt\Bigg]^{ \frac{1}{2}}
\cr
&\le C K_{ q}^{ \frac{3}{2}} A_q(1, 0)^{ \frac{3}{2}} r^{ -1} \Bigg[ \intl_{I_1} \| u(t)\|_{ L^q(B_1)}^3   dt\Bigg]^{ \frac{1}{2}}
\cr
&\le C r^{ -1}K_{ q}^{ \frac{3}{2}} A_q(1, 0)^{ 3}.
\label{4.27}
\end{align}
Furthermore, applying Poincar\'e's inequality, and employing Lemma\,\ref{lem4.6}, we find 
\begin{align}
&\intmw_{Q_r(z_0)}  | \nabla p_h \otimes \nabla p_h - (\widetilde {\nabla p_h \otimes \nabla p_h} )_{ B_r}|^{ \frac{3}{2}} dxdt
 \cr
 &\qquad \le C r^{-5 + \frac{3}{2}}  \int_{Q_r(z_0)}  | \nabla p_h|^{ \frac{3}{2}} | \nabla^2 p_h|^{ \frac{3}{2}} 
\le C r^{ - \frac{1}{2}}\intl_{Q_{ \frac{3}{4}}} | \nabla p_h| ^{3}
\cr
&\qquad \le C r^{ - \frac{1}{2}} A_q(1, 0)^3. 
\label{4.28}
\end{align}
By the aid of  \eqref{4.26a}, \eqref{4.27} and \eqref{4.28}  together with \eqref{4.25} 
we obtain  for all $ r_{ n+1} \le  r \le 1$
\[
\intl_{Q_r(z_0) }| u \otimes  u - (\widetilde {u \otimes  u} )_{ B_r(x_0)} |^{ \frac{3}{2}} \le C K_{ q}^3 A_q(1, 0)^3r^{ 4}.  
\]
Applying the Lemma\,\ref{lem4.8}, we find that for all  $ r_{ n+1} \le  r \le  r_2$
\begin{equation}
\intl_{Q_r(z_0)} | p_{ 2}- ({\tilde p} _2)_{ B_r(x_0)}|^{ \frac{3}{2}} \le C K_{ q}^3  A_q(1, 0)^3 r^{ 4}.  
\label{4.29}
\end{equation}

\hspace{0.5cm}
By $ \Psi _{ n+1}$ we denote the fundamental solution of the backward heat equation 
having its singularity at $ (x_0, t_0 + r_{ n+1}^2)$, more precisely,  
\[
\hspace*{-0.8cm}\Psi_{ n+1} (x,t) = \frac{c_0}{(r_{ n+1}^2- t+ t_0 )^{ \frac{3}{2}}} {\rm exp}  \Big\{- \frac{| x-x_0|^2}{(r^2_{ n+1}- t+ t_0)}\Big\}, \quad  (x, t) \in \R^{3}\times (-\infty, t_0+r_{ n+1}).  
\]
Taking a suitable cut off function $ \chi \in C^{\infty}( \R^{n})$ for $Q_{ r_4}(z_0)\subset  Q_{ r_3}(z_0)$, we may insert 
$ \phi = \Phi_{ n+1} = \Psi_{ n+1} \chi $ into the local energy inequality \eqref{4.24}  to get   for almost all $ t\in (t_0-r_3^2, t_0)$
\begin{align}
& \frac{1}{2}\intl_{B_{ r_3}(x_0)} \Phi_{ n+1} (t)| v(t)|^2  + \intl_{t_0-r_3^2}^{t}\intl_{B_{ r_3}(x_0)}\Phi _{ n+1} | \nabla v|^2  
\cr
&\quad \le \frac{1}{2}   \intl_{t_0-r_3^2}^{t}\intl_{B_{ r_3}(x_0)} | v|^2  
\Big(\frac{\partial \Phi_{ n+1}}{\partial t} + 
\Delta \Phi _{ n+1}\Big)
+ \frac{1}{2} \intl_{t_0-r_3^2}^{t}\intl_{B_{ r_3}(x_0)} | v|^2 v\cdot \nabla \Phi _{ n+1} 
\cr
& \qquad  - \frac{1}{2} \intl_{t_0-r_3^2}^{t}\intl_{B_{ r_3}(x_0)}  | v|^2 \nabla p_h\cdot \nabla \Phi_{ n+1}  + \intl_{t_0-r_3^2}^{t}\intl_{B_{ r_3}(x_0)} (v \otimes v - v \otimes \nabla p_h  : \nabla ^2 p_h) \Phi_{ n+1}
\cr
&\qquad  + \intl_{t_0-r_3^2}^{t}\intl_{B_{ r_3}(x_0)} p_1 v \cdot \nabla \Phi_{ n+1}+ \intl_{t_0-r_3^2}^{t}\intl_{B_{ r_3}(x_0)} p_2 v \cdot \nabla \Phi_{ n+1}. 
\label{4.26}
\end{align}
Arguing as in \cite{caf},  the above inequality yields  
\begin{align}
& \esssup_{t\in (t_0-r_{ n+1}^2, t_0)}\intmw_{B_{ r_{ n+1}}(x_0)} | v(t)|^2  + r_{ n+1}^{ -3}\intl_{Q_{ r_{ n+1} }(z_0)}  | \nabla v|^2  
\cr
&\qquad   \lesssim    \intl_{Q_{ r_3}(z_0)} | v|^2 
\Big|\frac{\partial \Phi_{ n+1}}{\partial t} + \Delta \Phi _{ n+1} \Big|
+ \intl_{Q_{ r_3}(z_0)} | v|^3 | \nabla \Phi _{ n+1}| 
+ \intl_{Q_{ r_3}(z_0)}  | v|^2 | \nabla p_h|| \nabla \Phi_{ n+1}|  
\cr
& \qquad\qquad   +\intl_{Q_{ r_3}(z_0)} | v|^2 | \nabla ^2 p_h| \Phi _{ n+1} 
+ \intl_{Q_{ r_3}(z_0)} | v|| \nabla p_h|  | \nabla ^2 p_h|\Phi _{ n+1} 
+\intl_{Q_{ r_3}(z_0)}  p_1  v\cdot  \nabla \Phi _{ n+1}
\cr
&\qquad \qquad \qquad \qquad +\intl_{Q_{ r_3}(z_0)}  p_2  v\cdot  \nabla \Phi_{ n+1}
\cr
&\qquad = I+ II+ III+ IV+ V+ VI+ VII. 
\label{4.31}
\end{align}

(i) Obviously, as $\Big|\frac{\partial \Phi_{ n+1}}{\partial t} + \Delta \Phi _{ n+1} \Big| \le C$ in $ Q_{ r_3}(z_0)$, and using \eqref{4.8},  we see that  
\[
I  \le C  \| v\|^2_{ L^3(Q_{ r_3}(z_0))} \le  CA_q(1, 0)^2. 
\]

(ii) As $ | \nabla \Phi _{ n+1}| \le Cr^{ -4}_k$ in $ Q_{ r_k}(z_0)  \setminus Q_{ r_{ k+1}}(z_0)$ for all $ k=1, \ldots, n$, observing 
(2.23)$ _{ k}$, and employing \eqref{4.26}, we get 
\begin{align*}
II &= \sum_{k=3}^{n} \intl_{Q_{ r_k}(z_0)  \setminus Q_{ r_{ k+1}}(z_0)} | v|^3  | \nabla \Phi _{ n+1}| + 
\intl_{Q_{ r_{ n+1}}(z_0) } | v|^3  | \nabla \Phi_{ n+1}| 
\\
 &\le C K_{ q}^3A_q(1, 0)^3 \sum_{k=2}^{n} r_k^{ -4} r_k^{ 5} \le C K_{ q}^3A_q(1, 0)^3. 
\end{align*}

(iii) Similarly as in (ii),  
\begin{align*}
III &= \sum_{k=3}^{n} \intl_{Q_{ r_k}(z_0)  \setminus Q_{ r_{ k+1}}(z_0)} | v|^2 | \nabla p_h |  | \nabla \Phi _{ n+1}| + 
\intl_{Q_{ r_{ n+1}}(z_0)  } | v|^2 | \nabla p_h|  | \nabla \Phi _{ n+1}| 
\\
 &\le C K_{ q}^2 A_q(1, 0)^3\sum_{k=1}^{n} r_k^{ -4}  r_k^{ \frac{13}{3}} \le C K_{ q}^3 A_q(1, 0)^3. 
\end{align*}

(iv) As $ \Phi _{ n+1} \le Cr^{ -3}_k$ in $ Q_{ r_k}(z_0)  \setminus Q_{ r_{ k+1}}(z_0)$ for all $ k=1, \ldots, n+1$ together with  (2.23)$ _k$ and \eqref{4.26}   we get 
\begin{align*}
IV &= \sum_{k=3}^{n} \intl_{Q_{ r_k}(z_0)  \setminus Q_{ r_{ k+1}}(z_0)} | v|^2 | \nabla^2 p_h |  \Phi_{ n+1} + 
\intl_{Q_{ r_{ n+1}}(z_0) } | v|^2 | \nabla^2 p_h|  | \Phi _{ n+1}| 
\\
 &\le C K_{ q}^2 A_q(1, 0)^3 \sum_{k=2}^{n} r_k^{ -3}  r_k^{ \frac{13}{3}}\le C K_{ q}^3A_q(1, 0)^3. 
\end{align*}

(v) Similarly as in (vi) we estimate  
\begin{align*}
V &= \sum_{k=3}^{n} \intl_{Q_{ r_k}(z_0)  \setminus Q_{ r_{ k+1}}(z_0)} | v| | \nabla p_h || \nabla^2 p_h |  \Phi_{ n+1}+ 
\intl_{Q_{ r_{ n+1}}(z_0)  } | v| | \nabla p_h || \nabla^2 p_h |  \Phi _{ n+1} 
\\
 &\le CK_{ q} A_q(1, 0)^3 \sum_{k=2}^{n} r_k^{ -3}  r_k^{ \frac{11}{3}} \le CK_{ q}^3 A_q(1, 0)^3. 
\end{align*}

(vi) To estimate $ VI$ we argue as in \cite{caf}. Let $ \chi _k$ denote cut off functions,  suitable for 
$ Q_{ r_{ k+1}}(z_0) \subset Q_{ r_{ k}}(z_0)$, 
$ k=2, \ldots, n+1$. Then 
\begin{align*}
VI &= \intl_{Q_{ r_{ 2}}(z_0)} p_1  v\cdot \nabla \Phi _{ n+1}
\\
&= \sum_{k=3}^{n} \intl_{Q_{ r_k}(z_0)  \setminus Q_{ r_{ k+2}}(z_0)} p_1 v\cdot \nabla (\Phi _{ n+1} (\chi _k - \chi _{ k+1})) 
\\
 &\qquad \qquad +\intl_{Q_{ r_{ 2}}(z_0)} p_1  v\cdot \nabla (\Phi _{ n+1} (1- \chi _2))   + \intl_{Q_{ r_{ 2}}(z_0)} p_1 v\cdot 
 \nabla (\Phi_{ n+1} \chi _{ n+1}) 
\\
&= \sum_{k=3}^{n} \intl_{Q_{ r_k}(z_0)  \setminus Q_{ r_{ k+2}}(z_0)} (p_1 - ({\tilde p} _1)_{ B_{ r_k}(x_0)}) v\cdot \nabla
 (\Phi_{ n+1} (\chi _k - \chi _{ k+1})) 
\\
 &\qquad \qquad + \intl_{Q_{ r_{ 2}}(z_0)} p_1  v\cdot \nabla (\Phi _{ n+1} (1- \chi _2))   
 \\
& \qquad \qquad + 
 \intl_{Q_{ r_{ n+1}}(z_0)} (p_1 - ({\tilde p} _1)_{ B_{ r_{ n+1}}(x_0)}) v\cdot \nabla (\Phi_{ n+1} \chi _{ n+1}).  
\end{align*}
As $ | \nabla (\Phi_{ n+1}(\chi_k - \chi_{ k+1}))|  \le C  r_k^{ -4}$ for $ k=1, \ldots, n$, applying Poincar\'e's inequality, 
using the fact that  $ p_1 $ is harmonic, together  with (2.23)$ _k$ and \eqref{4.8}  we see that 
\begin{align*}
 &\intl_{Q_{ r_k}(z_0)  \setminus Q_{ r_{ k+2}}(z_0)} (p_1 - ({\tilde p} _1)_{ B_{ r_k}(x_0)}) v\cdot \nabla (\Phi_{ n+1} (\chi _k - \chi _{ k+1})) 
 \\
  &\qquad \qquad \le CK_{ q} A_q(1, 0) r_k^{ -4} r_k^{ 5} \bigg(\intl_{Q_{ 1/2}} p_1^2 \bigg)^{ \frac{1}{2}}   
\\
&\qquad \qquad \le C  K_{ q}  (A_q(1, 0)^2 +A_q(1, 0)^{ \frac{3q}{2q-3}})r_k 
\\
&\qquad \qquad \le C  K_{ q} A_q(1, 0)^2r_k.   
\end{align*}
Summation from $ k=3 $ to $ n$ yields 
\[
 \sum_{k=3}^{n}\intl_{Q_{ r_k}(z_0)  \setminus Q_{ r_{ k+2}}(z_0)} (p_1 - ({\tilde p} _1)_{ B_{ r_k}(x_0)}) v\cdot \nabla (\Phi _{ n+1} (\chi _k - \chi _{ k+1}))  \lesssim  C  K_{ q}A(1, 0)^2.   
\]

Similarly, we find 
\begin{align*}
& \intl_{Q_{ r_{ 2}}(z_0)} p_1 v\cdot \nabla (\Phi _{ n+1} (1-\chi _2))- \intl_{Q_{ r_{ n+1}}(z_0)}
 (p_1 - ({\tilde p} _1)_{ B_{ r_{ n+1}}(z_0)}) v\cdot 
  \nabla (\Phi _{ n+1} \chi_{ n+1})) 
\\
&\quad \le (1+r_{ n+1})C  K_{ q} A_q(1, 0)^2.   
\end{align*}
Thus,
\[
VI \le C K_{ q} A_q(1, 0)^2.
\]

(vii) Finally, arguing as in (vi), and making   use of  \eqref{4.29}, we estimate 
\[
VII \le C K_{ q}^3 A_q^3(1, 0). 
\] 

Thus, inserting the estimates of $ I, II, III, IV, V, VI$ and $ VII$ into the right-hand side of  \eqref{4.31}, we get a constant $ C_q>0$ independently of $ n$,  such that  
\begin{align}
 &\esssup_{t\in (t_0-r_{ n+1}^2, t_0)}\intmw_{B_{ r_{ n+1}}(x_0)} | v(t)|^2  + r_{ n+1}^{ -3}\intl_{Q_{ r_{ n+1}}(z_0)}  | \nabla v|^2  
\cr
&\qquad  \le C_q \Big( K_{ q}^3 A_q(1, 0)^3 + K_{ q} A_q(1, 0)^2\Big)
= \Big(C_q K_{ q} A_q(1, 0)    + \frac{C_q}{K_{ q}}\Big) K_{ q}^2 A_q(1, 0)^2.
\label{4.32}
\end{align}
Note that $ C_q \rightarrow +\infty$ as $ q \rightarrow \frac{3}{2}$. 

\hspace{0.5cm}
On the other hand, using a standard interpolation argument along with \eqref{4.32}, we arrive at  
\begin{align}
\intmw_{Q_{ r_{ n+1}}(z_0)} | v|^3 &\le C_0 \bigg[\esssup_{t\in (t_0-r_{ n+1}^2, t_0)}\intmw_{B_{ r_{ n+1}}(x_0)} | v(t)|^2  + r_{ n+1}^{ -3}\intl_{Q_{ r_{ n+1}}(z_0)}  | \nabla v|^2\bigg] ^{ \frac{3}{2}} 
\cr
&\le   \Big(C_0C_q K_{ q} A_q(1, 0)   + \frac{C_0C_q}{K_{ q}}\Big)^{ \frac{3}{2}} K_{ q}^2A_q(1, 0)^2 
\label{4.33}
\end{align}
with an absolute constant $ C_0>1$. Note that neither $ C_q$ nor  $C_0 $ depend on the choice of $ K_{ q}$.  Thus we may set 
\[
K_{ q} = 2C_qC_0,\qquad  \var _{ q} = \frac{1}{4C_q^2C_0^2}. 
\]
Accordingly, if $ A_q(1, 0) \le \var _{ q}$,   \eqref{4.33}  implies  
\[
\intmw_{Q^{ n+1}(z_0)} | v|^3 \le K_{ q}^3 A_q(1, 0)^3.
\]
Whence, by  induction  the assertion of the proposition is true.  \hfill \Beweisende

\vspace{0.5cm}  
\hspace{0.5cm}
{\bf Proof of  Theorem\,\ref{thm4.5}}:   Proposition\,\ref{prop4.9} implies for every   Lebesgue point 
$ z_0=(x_0, t_0)\in Q_{ \frac{1}{2}}(0,0)$ of $ | v|^3$, after letting $ n \rightarrow +\infty$ in (2.23)$ _n$,  that the estimate  following holds true 
\begin{equation} 
| v(x_0,t_0)| \le K_{ q} A_q(1, 0). 
\label{4.34}
 \end{equation}
By using the triangular inequality and the mean value property of harmonic functions, we get from \eqref{4.34} 
 for almost all $ (x,t) \in Q_{ \frac{1}{2}}(0,0)$
 \begin{align*}
| u(x,t)| \le K_{ q} A_q(1, 0) + | \nabla p_h(x,t)| &\le K_{q} A_q(1, 0) + c \| u(t)\|_{L^q(B_1)}.
 \\
&\le  (K_{ q}+c) A_q(1,0). 
 \end{align*} 
This  leads to 
\begin{equation}
\| u\|_{ L^\infty(Q_{ \frac{1}{2}})} \le c (K_{ q} +1) \esssup_{t\in I_1} \| u(t)\|_{ L^q(B_1)} = c (K_{ q} +1) A_q(1,0). 
\label{4.36}
\end{equation}

Finally the assertion \eqref{4.7}    follows from  \eqref{4.36} respectively  by using a routine scaling argument.  \hfill \Beweisende

\section{Proof of Theorem\,\ref{thm1.1}}
\label{sec:-?}
\setcounter{secnum}{\value{section} \setcounter{equation}{0}
\renewcommand{\theequation}{\mbox{\arabic{secnum}.\arabic{equation}}}}

Let $ y_0\in \R^{3}$ be fixed. Let $ 0<r \le 1$ be arbitrarily chosen. By means of H\"older's inequality,  Sobolev's inequality and \eqref{6.3} (cf. Lemma\,\ref{lem6.1}) we get 
\begin{align}
r^{ -2}\intl_{B_r(y_0)} | U|^2 dy  &\le  C \intl_{B_1(y_0)} | U|^2+ | \nabla U|^2 dy  
\le  C\Big(\|  U\|_{ L^q(B_1(y_0))}^2+ \| \Omega \|_{ L^2(B_1(y_0))}^2\Big). 
\label{2.2}
\end{align} 
This implies for any $ 1 \le s_0 \le \,\sqrt{1+ 2a}$
\begin{equation}
r^{ -2}\intl_{B_r( \frac{y_0}{s_0})} | U|^2 dy \le C \sup_{ 1 \le s \le \,\sqrt{1+2a}} 
\Big(\|  U\|_{ L^q(B_1(\frac{y_0}{s}))}^2+ \| \Omega \|_{ L^2(B_1(\frac{y_0}{s}))}^2\Big) =: C\Psi (y_0). 
\label{2.2}
\end{equation}

\hspace{0.5cm}
Recalling the definition of $ u$, and  using \eqref{2.2},  we get for almost all $ t \in \Big(- \frac{1}{2a}-1 , - \frac{1}{2a}\Big)$
\begin{align*}
\intl_{B_r(y_0)} | u(x,t)|^2 dx &= \frac{1}{ -2at} \intl_{B_r(y_0)} \Big| U \Big(\frac{x}{ \,\sqrt{-2at}}\Big)\Big|^2 dx
\\
&= \,\sqrt{-2at} \intl_{B_{ \frac{r}{ \,\sqrt{-2at}}}( \frac{y_0}{ \,\sqrt{-2at}})} |  U|^2 dy
\\
& \le \,\sqrt{1+ 2a}\intl_{B_{r}( \frac{y_0}{ \,\sqrt{-2at}})} |  U|^2 dy \le  C_0r^2 
\Psi (y_0), 
\end{align*}
Accordingly, setting $  z_0 = \Big(y_0, - \frac{1}{2a}\Big)$, the above estimate becomes 
\begin{equation}
A_2(r, z_0)^2 \le C_0 r \Psi (y_0).
\label{2.2a}
\end{equation}

We take 
\[
r= \min \bigg\{1,    \frac{\var _2^2}{C_0\Psi  (y_0)}\bigg\}.
\]
where $ \var _2$ denotes the constant in Theorem\,\ref{thm4.5} for the case $ q=2$. By the 
choice of $ r$ we infer from \eqref{2.2a}  
\begin{equation}
A_2(r, z_0) \le \var _2. 
\label{2.3}
\end{equation}
Accordingly, Theorem\,\ref{thm4.5} together with \eqref{2.2a}  yields 
\begin{align}
| U(y_0)| &\le \sup_{ Q_{ \frac{r}{2}}(z_0)} | u| \le C_2  r^{ -1} A_2(r, z_0) 
\le C_2 C_0^{ \frac{1}{2}}r^{ - \frac{1}{2}}  \Psi  (y_0)^{ \frac{1}{2}}
\cr
&\le C_2  (\var _2^{ -1} \Psi  (y_0) + C_0^{ \frac{1}{2}} \Psi  (y_0)^{ \frac{1}{2}}). 
\label{2.4}
\end{align}
By the  assumption \eqref{1.6}   having $ \Psi  (y_0) = {\rm o}(| y_0|)$, it follows from \eqref{2.4} that 
$ U$ has sublinear growth.  Furthermore, appealing to Lemma\,\ref{lemA.2}, we see that $ | P(y)|= O(| y|^{ \frac{9}{2}})$ which allows to apply  the maximum principle for the energy function $ \Pi = \frac{| U|^2}{2}+ P + ay\cdot U $.   Hence, $ \Pi = \Pi _0 =\const$ (cf. Lemma\,5.1\cite{tsa}). Thus,  applying the formula
$$
-\Delta \Pi +(U+ay)\cdot \nabla \Pi =-|\Omega|^2, \quad \Omega=\nabla \times U,
$$
 we see that $ \Omega=0$. Combining this with the condition $\nabla \cdot U=0$, we find that
each component $ U_i$, $ i=1,2,3$, is  harmonic. Since $ U$ has sublinear   growth  at infinity, 
we get $ U=\const$.

 \hfill \Beweisende 

\section{Proof of Theorem\,\ref{thm1.2}}
\label{sec:-7}
\setcounter{secnum}{\value{section} \setcounter{equation}{0}
\renewcommand{\theequation}{\mbox{\arabic{secnum}.\arabic{equation}}}}

  Let $ y_0\in \R^{3}$ be fixed. Let $ \alpha >0$ be chosen so that \eqref{1.8} holds true.  
Recalling the definition of $ u$,  we get for any $ 0< r \le 1$ and for  almost all $ t \in \Big(- \frac{1}{2a}-1 , - \frac{1}{2a}\Big)$
\begin{align}
&\intl_{B_r(y_0)} | u(x,t)|^q dx 
\cr
&\quad = \frac{1}{ (\,\sqrt{-2at})^{ q}} \intl_{B_r(y_0)} \Big| U \Big(\frac{x}{ \,\sqrt{-2at}}\Big)\Big|^q dx
\cr
&\quad  = (\,\sqrt{-2at})^{ 3-q}\intl_{B_{ \frac{r}{ \,\sqrt{-2at}}}( \frac{y_0}{ \,\sqrt{-2at}})} |  U|^q dy
\cr
&\quad  \le \alpha ^q\mes( B_1) r^3 +  (1+ 2a)^{ \frac{3-q}{2}} \intl_{B_{r}(\frac{y_0}{ \,\sqrt{-2at})})  \cap \{ | U|>\alpha\} }   |  U|^q dy 
\cr
&\quad = C_0 r^3 + {\tilde \Phi}(y_0),
\label{7.0}
\end{align}
where we set $ {\tilde \Phi } (y_0)=(1+ 2a)^{ \frac{3-q}{2}}  \sup_{ 1 \le s \le \,\sqrt{1+2a}}  \| U\|^q_{ L^q(B_1 (\frac{y_0}{s}))}
$ and $ C_0= \alpha ^q\mes (B_1)$. 
Setting $ z_0= \Big(y_0, - \frac{1}{2a}\Big)$ from the inequality above,  we deduce that 
\begin{equation}
A_q(r, z_0)^q = \sup_{t\in  (- \frac{1}{2a}-r^2 , - \frac{1}{2a}) } \| u(t)\|^q_{ L^q(B_1(y_0))} \le C_0 r^q + r^{ q-3}{\tilde \Phi}(y_0).
\label{7.1}
\end{equation}

\hspace{0.5cm}
Without loss of generality we may assume that $ C_0>1$. We now take $ r$ such that 
\begin{equation}
C_0 r^q = \frac{\var _q^q}{2},
\label{7.4}
\end{equation}
where $ \var _q$ denotes the positive number in Theorem\,\ref{thm4.5}. 
On the other hand,  our assumption \eqref{1.7} yields ${\tilde \Phi}(y_0) \to 0$ as $ | y_0| \rightarrow +\infty$.  
Therefore, we may chose $ R>0$ such that for all $ y_0 \in \R^{3}  \setminus B_R$ 
\[
r^{ q-3}{\tilde \Phi}(y_0) \le \frac{\var _q^q}{2},
\]
 Accordingly, Theorem\,\ref{thm4.5} implies for all $ y_0 \in \R^{3}  \setminus B_R$
\begin{equation}
| U(y_0)| \le \sup_{ Q_{ \frac{r}{2}}(z_0)} | u| \le C_q  r^{ -1}A_q(r, z_0) \le C_q r^{ -1} \var _q. 
\label{7.2}
\end{equation}
Therefore $ U$ is bounded. According  to Tsai's result (cf. Lemma\,5.1\cite{tsa}), we conclude that $ U=\const$. This completes the 
proof of the theorem.  \hfill \Beweisende 

\section{Proof of Theorem\,\ref{thm1.3}}
\label{sec:-8}
As in \cite{cha} we consider the self-similar transform of the solution $(u,p)$ of (\ref{1.1})  into $(U,\Pi)$ by
$$
u(x,t)=\frac{1}{\sqrt{2a(t_*-t)}} V(y,s), \quad p(x,t)=\frac{1}{2a(t_*-t)} \Pi(y,s),
$$
where
$$
y=\frac{x-x_*}{\sqrt{2a(t_*-t)}} , \quad s=\frac12 \log \left(\frac{t_*}{t_* -t}\right).
$$
Then, the system  (\ref{1.1}) is transformed into a system for $(V,\Pi)\in C^2 (\Bbb R^3 \times (0, \infty))$
\begin{equation}
V_s-\Delta V  + (V\cdot \nabla ) V + a y\cdot \nabla V + a V = - \nabla \Pi,\qquad  \nabla \cdot  V =0\quad  \text{ in}
\quad  \R^{3}. 
\label{7.5}
\end{equation}
The condition (\ref{1.9}) is transformed into 
\begin{equation}
\lim_{s\to \infty} \|V(\cdot ,s)-U\|_{L^q \left(B_{\frac{r}{\sqrt{2a}}} (0)\right)}=0 \quad \forall r >0.
\label{7.6}
\end{equation}
From the argument of Proof of   \cite[Theorem\,1.2]{cha})  one can show from (\ref{7.6}) that
$U$ is a solution of (\ref{1.4}) for a scalar function $P$. We include this part here for reader's convenience.
We choose  $\xi \in C^1_c(0,1)$ with $\int_0 ^1 \xi (s)ds =1$, and $\varphi \in C_c ^1 (\Bbb R^3)$ with $\nabla \cdot \phi =0$.
Then, multiplying (\ref{7.5}) by $\xi (s-n)\varphi(y)$ and integrating it over $\Bbb R^3 \times [n, n+1]$, then after integration by part we obtain 
\bq
\lefteqn{\int_0 ^1\int_{\Bbb R^3} \xi '(s) \varphi (y) \cdot V(y, s+n) dyds +\frac{a}{2}\int_0 ^1 \int_{\Bbb R^3} \xi (s)  V(y, s+n) \cdot \varphi (y) dy ds}\nonumber \\
&&=-\int_0 ^1 \int_{\Bbb R^3} \xi (s) (V\otimes V)(y, s+n) : \nabla  \varphi (y) dyds \nonumber \\
&& \quad -\int_0 ^1 \int_{\Bbb R^3} \xi (s) \left\{a V(y, s+n)\cdot  (y\cdot \nabla) \varphi  -V(y,s+n)\cdot \Delta \varphi \right\} dyds.
\label{7.7}
\eq
Since (\ref{7.6}) implies that $V(\cdot, s+n) \to U$ in $L^2_{\rm{loc}} (\Bbb R^3)$ as $n\to \infty$, passing $n\to \infty$ in \eqref{7.7}, using the fact $\int_0 ^1 \xi (s)ds =1$, one has
\bq
\frac{a}{2}\int_{\Bbb R^3}  \varphi (y) \cdot U dy &=&- \int_{\Bbb R^3}  (U\otimes U)(y) : \nabla  \varphi (y) dy\nonumber \\
 &&- \int_{\Bbb R^3} \left\{a U(y)\cdot  (y\cdot \nabla) \varphi  -U(y)\cdot \Delta \varphi \right\} dy,
\label{7.8}
\eq
which shows that $U\in L^2_{\rm{loc}} (\Bbb R^3)$ is a weak solution of (\ref{1.4}) for some scalar function $P=P(y)$. By a standard regularity theory $(U,P)$ is a smooth solution of 
(\ref{1.4}). Now, applying Corollary 1.4 one can conclude that 
$U=0$, and the condition (\ref{1.9}) reduces to 
\begin{equation}
\label{7.9}
\lim_{t\to t_*}(t_*-t)^{\frac{q-3}{2q} } \sup_{t<\tau <t_*}\left \|u (\cdot , \tau) \right\|_{L^q (B_{r\sqrt{t_*-t} } (x_*))} =0
\end{equation}
for each $r>0$.  Setting $r=1$, $\rho=\sqrt{t_*-t}$,  we find that
$$
\lim_{\rho\to 0} \left\{ \rho^{\frac{q-3}{q}} \sup_{t_* -\rho^2 <\tau < t_* }\|u(\cdot , \tau )\|_{L^q (B(x_*, \rho ))} \right\} = \lim_{\rho\to 0} A_q(z_{ \ast}, \rho )=0,
$$
where $ z_*=(x_*, t_*)$. 
Thanks to Theorem\,2.5  we find that $z_*$ is a regular point
(cf. also  the regularity criterion  due to Gustafson, Kang and Tsai \cite[Theorem 1.1]{gus} ).
  \hfill \Beweisende 

 \hspace{0.5cm}
$$\mbox{\bf Acknowledgements}$$
Chae was partially supported by NRF grants 2016R1A2B3011647, while Wolf has been supported 
supported by the German Research Foundation (DFG) through the project WO1988/1-1; 612414.  

\appendix

\section{Gradient estimates and pressure estimate}
\label{sec:-A}
\setcounter{secnum}{\value{section} \setcounter{equation}{0}
\renewcommand{\theequation}{\mbox{A.\arabic{equation}}}}

\begin{lem}
\label{lem6.1}
Let $ U\in L^q_{ loc}(\R^{3})$, $ 0<q<+\infty$, with $ \nabla\cdot  U=0$. 
Let $ \Omega = \nabla \times U\in L ^2_{ loc}(\R^{3})$.  Then $ \nabla u\in L^2_{ loc}(\R^{3})$, and there holds 
for all $ y_0 \in \R^{3}$
\begin{equation}
\| U\|_{ L^2(B_1(y_0))}+ \| \nabla U\|_{ L^2(B_1(y_0))}  \lesssim  
 \| U\|_{ L^q(B_2(y_0))} + \|\Omega \|_{ L^2(B_2(y_0))}, 
\label{6.3}
\end{equation}  
where the hidden constant depends only on $ q$. 
\end{lem}

{\bf Proof}:  By means of  a standard mollifying argument it suffice to verify the estimate \eqref{6.3} for smooth $ U$.  

\hspace{0.5cm}
Let $ \alpha = \frac{6-q}{2q}$. For $ y_0\in \R^{3}$ let $ \zeta \in C^{\infty}_{\rm c}(B_2(y_0))$ denote a suitable cut off function for $ B_1(y_0) \subset B_2(y_0)$.  
Using integration by parts, we find 
\begin{align*}
&\intl_{B_2(y_0)} | \nabla U|^2 \zeta ^{ 2\alpha } dx = \intl_{B_2(y_0)} \nabla U:\nabla U \zeta ^{ 2\alpha } dx 
\\
&\qquad = -\intl_{B_2(y_0)}  U\cdot \Delta  U \zeta ^{ 2\alpha } dx - \alpha \intl_{B_2(y_0)}  \nabla | U|^2\cdot  \zeta ^{ 2\alpha -1} \nabla \zeta  dx
\\
&\qquad = \intl_{B_2(y_0)}  U\cdot \nabla \times  \Omega  \zeta ^{ 2\alpha } dx +\alpha  \intl_{B_2(y_0)}  | U|^2 \zeta ^{ 2\alpha -1}  \Delta \zeta  dx
\\
& \qquad \qquad \qquad +\alpha  (2\alpha -1)\intl_{B}  | U|^2\zeta ^{ 2\alpha -2}| \nabla \zeta |^2 dx. 
\\
&\qquad = \intl_{B_2(y_0)}  | \Omega |^2\zeta ^{ 2\alpha } dx- 2\alpha \intl_{B_2(y_0)}  
 \Omega\cdot \zeta ^{ 2\alpha -1} U\times \nabla \zeta  dx 
\\
& \qquad \qquad \qquad +\alpha  \intl_{B_2(y_0)}  | U|^2 \zeta ^{ 2\alpha -1}  \Delta \zeta  dx+\alpha  (2\alpha -1)\intl_{B_2(y_0)}  | U|^2\zeta ^{ 2\alpha -2}| \nabla \zeta |^2 dx. 
\end{align*}
Applying Young's inequality, we get 
\begin{align}
&\intl_{B_2(y_0)} | \nabla U|^2 \zeta ^{ 2\alpha } dx   \lesssim   \intl_{B_2(y_0)} | \Omega |^2 \zeta ^{ 2\alpha } dx + \intl_{B_2(y_0)} | U |^2 \zeta ^{ 2\alpha -2} dx.
\label{6.1}
\end{align}
Applying  H\"older's inequality along with Sobolev's inequality we estimate 
\begin{align*}
&\| U \zeta^{ \alpha -1} \|_{ L^2(B_2(y_0))}^2 
\\
&\quad \le \| U\|_{ L^q(B_2(y_0))}^{ \frac{4q}{6-q}} 
\| U\zeta ^{ (\alpha -1)\frac{6-q}{6-3q}}\|_{ L^6(B_2(y_0))}^{ \frac{12-6q}{6-q}} 
= \| U\|_{ L^q(B_2(y_0))}^{ \frac{4q}{6-q}} 
\| U \zeta ^{\alpha }\|_{ L^6(B_2(y_0))}^{ \frac{12-6q}{6-q}} 
\\
& \quad  \lesssim   \| U\|_{ L^q(B_2(y_0))}^{ \frac{4q}{6-q}} 
\| \nabla U \zeta ^{\alpha }\|_{ L^2(B_2(y_0))}^{ \frac{12-6q}{6-q}} + \| U\|_{ L^q(B_2(y_0))}^{ \frac{4q}{6-q}} 
\| U \zeta ^{\alpha -1}\|_{ L^2(B_2(y_0))}^{ \frac{12-6q}{6-q}}dx.
\end{align*}
Applying Young's inequality, we are led to 
\begin{equation}
\| U \zeta^{ \alpha -1} \|_{ L^2(B_2(y_0))}^2  \lesssim    \| U\|_{ L^q(B_2(y_0))}^{ \frac{4q}{6-q}} 
\| \nabla U \zeta ^{\alpha }\|_{ L^2(B_2(y_0))}^{ \frac{12-6q}{6-q}} + \| U\|_{ L^q(B_2(y_0))}^{2}.  
\label{6.2}
\end{equation}
Combining \eqref{6.1} and \eqref{6.2}, and using once more Young's inequality, we obtain \eqref{6.3}.  

\hfill \Beweisende

\begin{lem}
\label{lemA.2}
Let $ (U,P)$ be a smooth solution to \eqref{1.4}. Assume that $ | U(y)|=O(| y|)$, and $ \| \nabla U\|_{ L^2(B_1(y))}= 
O(| y|^{ \frac{1}{2}})$ as $|y|\to \infty$. Then, we have
\begin{align}
\| \nabla^k U\|_{ L^2(B_1(y))}&= O(| y|^{ \frac{2k-1}{2}}),
\label{6.6a}\\
 | \nabla^k U(y)|&= O(| y|^{ \frac{2k+3}{2}}),
   \label{6.6c} \\
| P(y)|&= O(| y|^{ \frac{9}{2}})
\label{6.6b}
\end{align} 
as $|y|\to \infty$.
\end{lem}

{\bf Proof}: Appying $ \nabla \times $ to both sides of \eqref{1.4}, we obtain  
\begin{equation}
- \Delta  \Omega + \nabla \times ((U\cdot \nabla) U) + a y\cdot \nabla \Omega + 2a \Omega = 0. 
\label{6.3a}
\end{equation}
Let $ y_0\in \R^{3}$, and let $ \zeta \in C^{\infty}_{\rm c}(B_2(y_0))$ denote a suitable  cut off function for 
$ B_1(y_0 ) \subset  B_2(y_0)$. Multiplying \eqref{6.3a} by $ \Omega  \zeta ^2$, and applying integration by parts,  
using the formula $U\cdot \nabla U =\frac12 \nabla |U|^2 +\Omega \times U$, we get 
\begin{align}
&\intl_{B_2(y_0)} | \nabla \Omega |^2 \zeta ^2 + \frac{a}{2}| \Omega |^2 \zeta ^2 dx 
\cr
&\quad = \frac{1}{2}   \intl_{B_2(y_0)} |  \Omega |^2 
(\Delta \zeta ^2  + a y\cdot \nabla \zeta ^2 )dx
- \intl_{B_2(y_0)} \Omega \times U \cdot  \nabla \times (\Omega   \zeta ^2)   dx. 
\label{6.4}
\end{align}
Using Cauchy-Schwarz's inequality and Young's inequality, we find 
\begin{align}
\intl_{B_2(y_0)} | \nabla \Omega |^2 \zeta ^2 + 2a| \Omega |^2 \zeta ^2 dx \le C (1+ | y_0|+ | y_0|^3). 
\label{6.5}
\end{align}
With the aid of  Lemma\,\ref{lem6.1} together with Sobolev's inequality we infer from \eqref{6.5}
\begin{align}
\| \nabla ^2 U\|_{ L^2(B_1(y_0))}+ \|\nabla  U\|_{ L^6(B_1(y_0))}  \le  C (1+ | y_0|^{ \frac{3}{2}}). 
\label{6.6}
\end{align}  
Differentiating  \eqref{6.3},  arguing as above together, and applying  an inductive argument, we see that   
\begin{align}
\| \nabla ^k U\|_{ L^2(B_1(y_0))}+ \|\nabla^{ k-1}  U\|_{ L^6(B_1(y_0))}  \le  C (1+ | y_0|^{ \frac{2k-1}{2}}). 
\label{6.7}
\end{align}  
Whence, \eqref{6.6a}.  Employing Sobolev's embedding theorem, we get $ H^{ k+2}(B_1(y_0)) \hookrightarrow L^\infty(B_1(y_0))$.  Thus \eqref{6.6c} is an immediate consequence of \eqref{6.6a} by using  Sobolev's inequality.  

\hspace{0.5cm}
Observing  \eqref{1.4}, and making use of  \eqref{6.6c}, we obtain 
\begin{equation}
| \nabla P(y) | = O(| y|^{ \frac{7}{2}}).
\label{6.10}
\end{equation}  
This immediately implies \eqref{6.6b}.  \hfill \Beweisende

\bibliographystyle{siam}

\end{document}